\documentclass[12pt]{iopart}
\usepackage{mathrsfs, amsthm, amssymb, harvard}

\newtheorem{theorem}{Theorem}[section]
\newtheorem{lemma}[theorem]{Lemma}
\newtheorem{proposition}[theorem]{Proposition}

\newtheorem{corollary}[theorem]{Corollary}

\def\be#1{\begin{equation}\label{#1}}
\def\ee{\end{equation}}
\def\bea{\begin{eqnarray*}}
\def\eea{\end{eqnarray*}}
\def\e{{\rm e}}
\def\d{{\rm d}}
\def\uu{u}\def\ux{x}\def\uy{y}
\def\half{\frac{1}{2}}
\def\ddt{\frac{\d}{\d t}}
\def\lap{\Delta}
\def\grad{\nabla}
\def\Re{{\mathbb R}}

\begin{document}

\title[Singularities in the 3D NSE]{Decay of weak solutions and the singular set of the three-dimensional Navier--Stokes equations}

\author{James C. Robinson}

\address{Mathematics Institute, University of Warwick, Coventry CV4 7AL, UK}
\ead{j.c.robinson@warwick.ac.uk}

\author{Witold Sadowski}
\address{$^1$ Faculty of Mathematics, Informatics and Mechanics,
       Warsaw University,
       Banacha 2 02-097, Warszawa, Poland.}
\address{$^2$ Institute of Mathematics of the Polish Academy of Sciences,
       ul. \'Sniadeckich 8 00-956, Warszawa, Poland.} \ead{witeks@hydra.mimuw.edu.pl}

\begin{abstract}
We consider the behaviour of weak solutions of the unforced
three-dimensional Navier-Stokes equations, under the assumption that
the initial condition has finite energy ($\|u\|^2=\int|u|^2$) but
infinite enstrophy ($\|Du\|^2=\int|{\rm curl}\,u|^2$). We show that
this has to be reflected in the solution for small times, so that in
particular $\|Du(t)\|\rightarrow+\infty$ as $t\rightarrow0$. We also
give some limitations on this `backwards blowup', and give an
elementary proof that the upper box-counting dimension of the set of
singular times can be no larger than one half. Although similar in flavour, this final result neither implies nor is implied by Scheffer's result that the $1/2$-dimensional Hausdorff measure of the singular set is zero.
\end{abstract}

\submitto{\NL}
\maketitle

\section*{Introduction}
The existence or otherwise of singularities in the three-dimensional
Navier--Stokes equations is acknowledged to be one of the
outstanding open problems in the theory of partial differential
equations, and so there has been much research concerned with the
time before the formation of a possible singularity. In this paper
we consider the behaviour of solutions after the occurrence of a
singularity, or equivalently (but less hypothetically) the time
evolution of a weak solution with a suitably irregular initial
condition.

In order to describe our results in more detail we need some simple
terminology and notation. Throughout the paper we consider the
evolution of an incompressible fluid in a periodic box
$\Omega=[0,L]^3$ in the absence of external forces: the evolution of
the velocity field $u(x)$ is modelled by the three-dimensional
Navier--Stokes equations
$$
\frac{\partial u}{\partial t}-\nu\Delta u+(u\cdot\nabla)u+\nabla
p=0\qquad\nabla\cdot u=0.
$$
For convenience sake we also assume that the total momentum is zero, $\int_\Omega u(x)\,\d x=0$, a situation that is preserved under the evolution of $u$.

We denote the kinetic energy of $u(x)$ by
$$
\|u\|^2:=\int_\Omega |u(x)|^2\,\d x
$$
and the `enstrophy' of $u(x)$, which is the square integral of the
vorticity $\omega(x)={\rm curl}\ u(x)$, by
$$
\|Du\|^2:=\int_\Omega|{\rm curl}\,u(x)|^2\,\d x.
$$
We denote the space of all divergence-free flow fields with finite
kinetic energy by $H$, and the space of all such flow fields with
finite enstrophy by $V$. (More strictly, if $\mathscr V$ is the
space of all $C^\infty$ $\Omega$-periodic functions that are
divergence free and have zero average over $\Omega$, then $H$ is the closure of $\mathscr V$ in the
$L^2$ norm, and $V$ is the closure of $\mathscr V$ in the $H^1$
norm.)

Leray showed in 1934 that given an initial condition with finite
energy, there is a weak solution of the equations that exists for
all time and whose kinetic energy remains finite. Indeed, we have
$$
\|u(t)\|^2+\int_0^t\|Du(s)\|^2\,\d s\le\|u_0\|^2
$$
for all $t\ge0$, and the solution is weakly continuous into $H$,
i.e.~$u(t)\rightharpoonup u(s)$ as $t\rightarrow s$. However, the
uniqueness of such a solution is not guaranteed.

On the other hand, while it is know that if the enstrophy
$\|Du(t)\|^2$ remains finite then the solution is unique, an initial
condition $u_0$ with $\|Du_0\|^2<\infty$ is only known to give rise
to a solution whose enstrophy remains finite for a limited time.
[For an elementary discussion of these issues see Robinson (2006);
for a more advanced but accessible treatment see Doering \& Gibbon;
and for a more technical mathematical discussion see Constantin \&
Foias (1988) or Temam (1977).]

The possibility of singular behaviour in the Navier-Stokes equations
was first pointed out by Leray (1934) who suspected the existence of
self-similar singular solutions. He also showed (using a simple
argument which we will review later) that any solution that blows up
at a time $t=T$ must be bounded below by
$$
\|Du(t)\|^2\ge\frac{c}{\sqrt{T-t}}.
$$

Still little is known of the regularity (or otherwise) of weak
solutions of the Navier-Stokes equations, but several papers have
imposed some restrictions on the hypothetical blow up. The most
famous is probably Scheffer's result stating that the Hausdorff
measure of turbulent times in the Navier-Stokes equations is less
than one half. More sophisticated is the result of Caffarelli, Kohn,
\& Nirenberg (1982), which treats the singular set in space-time
(i.e.~as a subset of $\Omega\times[0,\infty)$). Scheffer's result
can be obtained as a consequence of the result of Caffarelli et al.,
see Section 7 in Robinson (2006).

Here we give a simple proof of a result reminiscent of Scheffer's:
we show that the upper box-counting dimension of the singular set is
no greater than $1/2$. This result neither implies, nor is implied
by, Scheffer's.

\section{Solutions with irregular initial data}

Given initial data $u_0\in V$, it is well known that the solution
remains regular (i.e.~$\|Du(t)\|$ remains finite) at least for some
short time interval. Indeed, standard estimates can be used to show
that for $t>s$
  \be{upper}
  \|Du(t)\|^2\le\frac{\|Du(s)\|^2}{\sqrt{1-c(t-s)\|Du(s)\|^4}}
  \ee
  while the right-hand side remains finite, i.e.~while
  \be{regulartime}
  c(t-s)\|Du(s)\|^4<1.
\ee
In particular, if $s=0$ then this shows that the solution remains
regular at for $t\in[0,T)$, where $T=c^{-1}\|Du_0\|^{-4}$.

 Rearranging (\ref{upper}) to obtain a lower bound for $s<t$ gives
  \be{lower}
  \|Du(s)\|^2\ge\frac{\|Du(t)\|^2}{\sqrt{1+c(t-s)\|Du(t)\|^4}}.
  \ee
While one might be concerned as to the range of validity of this
bound (given the conditions on $t$ and $s$ in (\ref{regulartime}),
it does in fact hold for all $s\le t$. Indeed, if for some $s<t$ one
has
    $$
   \|Du(s)\|^2<\frac{\|Du(t)\|^2}{\sqrt{1+c(t-s)\|Du(t)\|^4}},
$$
then the right-hand side of (\ref{upper}) remains finite on $(s,t]$
and one can infer that $\|Du(t)\|^2<\|Du(t)\|^2$, a contradiction.

From (\ref{lower}) one can easily deduce (as did Leray) that if a
solution does blow up as $t\rightarrow T$, i.e.~if
$\|Du(t)\|\rightarrow\infty$ as $t\rightarrow T$, then for $s<T$ one
must have the lower bound
\be{Lerays}
\|Du(s)\|^2\ge\frac{1}{\sqrt{c(T-s)}}\qquad\mbox{for all}\quad 0\le
s<T.
\ee

Given $u_0\in V$, we know that even if there is a loss of regularity
at some time $T$, there is a weak solution that exists for all
positive times. It is therefore reasonable to ask about the
behaviour of this solution {\it after} the putative blowup time.
Even the equations are in fact `regular' - i.e.~if any initial
condition in $V$ gives rise to a regular solutions that is valid for
all time - one can still ask the same question about solutions whose
initial condition is contained in $H$ and not in $V$.

Our first simple result shows that such solutions retain a trace of
this `irregularity', i.e.~it is not possible that $\|Du_0\|=\infty$
but $\|Du(t)\|$ is uniformly bounded on an interval of the form
$(0,\epsilon)$ for some $\epsilon>0$.

\begin{lemma}
  Suppose that $u_0\in H\setminus V$, and that $u(t)$ is a weak solution
  with $u(0)=u_0$. Then $\|Du(t)\|\rightarrow+\infty$ as $t\rightarrow0$.
\end{lemma}

\begin{proof}
Suppose not. Then for some $M>0$ we have a sequence $t_n$ with
$\|Du(t_n)\|\le M$. Since $u$ is weakly continuous into $H$, we know
that $u(t_n)\rightharpoonup u_0$ in $H$. However, since
$\|Du(t_n)\|\le M$ it follows from the Alaoglu Theorem that
$u(t_n)\rightharpoonup v_0$ in $V$, for some $v_0$ with $\|Dv_0\|\le
M$. Since weak convergence in $V$ implies weak convergence in $H$,
and weak limits are unique, we must in fact have $u_0=v_0$, and so
$\|Du_0\|\le M$, a contradiction.
\end{proof}

\section{Decay of solutions}

Given the result of this lemma, it becomes natural to look for
possible upper or lower bounds on the rate of decay of the $V$ norm
for solutions with irregular initial data, and we now turn to some
results in this direction.

All of our results will be consequences of the lower bound in
(\ref{lower}), which is more convenient to work with when written in
the form
  \be{lower2}
  \|Du(s)\|^2\ge\frac{1}{\sqrt{\|Du(t)\|^{-4}+c(t-s)}}.
  \ee
We note here that this lower bound remains valid even if $u$ is not regular on the interval $[s,t]$, since $\|Du(r)\|=+\infty$ for $r\in(s,t)$ implies that
$$
\|Du(s)\|^2\ge\frac{1}{\sqrt{c(r-s)}}\ge\frac{1}{\sqrt{\|Du(t)\|^{-4}+c(t-s)}}.
$$

\subsection{Solutions in $L^4(0,T;V)$}

An immediate consequence of the inequality in (\ref{lower2}) is an
upper bound on the $V$ norm of solutions known to be in
$L^4(0,T;V)$. In order for
$$
\int_0^T\|Du(s)\|^4\,\d s
$$
to be finite, one would expect $\|Du(t)\|^4$ blows up more slowly
than $t^{-1}$ as $t\rightarrow0$. In particular, it is natural to
imagine that $t\|Du(t)\|^4$ remains bounded as $t\rightarrow0$, and
this is shown rigorously in the following result:

\begin{proposition}
Suppose that $u\in L^4(0,T;V)$. Then there exists a constant $c$
depending on $\|u\|_{L^4(0,T;V)}$ such that
\be{L42upper}
t^{1/2}\|Du(t)\|^2\le c.
\ee
\end{proposition}

\begin{proof}
  From (\ref{lower}) we have
$$
\|Du(s)\|^4\ge\frac{\|Du(t)\|^4}{1+c(t-s)\|Du(t)\|^4},
$$
and so
\bea
\int_0^t\|Du(s)\|^4\,\d
s&\ge&\left[-\frac{1}{c}\log(1+c(t-s)\|Du(t)\|^4)\right]_0^t\\
&=&\frac{1}{c}\log(1+ct\|Du(t)\|^4)
\eea
from which (\ref{L42upper}) follows.
\end{proof}

It is natural to conjecture that this decay also holds for all weak
solutions (this is the case for the two-dimensional Navier--Stokes
equations), but in general all the information we currently have
available is that $u\in L^2(0,T;V)$, and so a more natural
conjecture along similar lines is that $t\|Du(t)\|^2\le c$. However,
a similar approach, but starting from (\ref{lower2}) yields much
weaker results. Indeed, if $\|Du(t_n)\|=+\infty$ with
$t_n\downarrow0$ then one can only obtain
$$
\sum_{n=1}^\infty(t_n-t_{n+1})^{1/2}<\infty,
$$
an observation originally due to Leray.

\subsection{Box-counting dimension and decay of solutions}

Given a set $X$, its (upper) box-counting dimension is defined as
follows. Let $N(X,\epsilon)$ denote the largest number of disjoint
balls of radius $\epsilon$ with centres in $X$. Then
$$
d_f(X)=\limsup_{\epsilon\rightarrow0}\frac{\log
N(X,\epsilon)}{-\log\epsilon},
$$
i.e.~essentially $N(X,\epsilon)\sim\epsilon^{-d_f(X)}$. [Usually one
would take $N(X,\epsilon)$ to be the smallest number of balls of
radius $\epsilon$ whose union covers $X$ (see Definitions 3.1 in
Falconer, 1990), but the two definitions are equivalent and that
given here is more convenient for what follows.]

We now give an indication of how the box-counting dimension occurs
naturally when considering rates of decay of solutions. Suppose that
\be{rate}
\|Du(t)\|^2\ge ct^{-\alpha}.
\ee
Then, since we know that $\int_0^T\|Du(s)\|^2\,\d s\le\|u_0\|^2$, it
must be the case that $\alpha<1$. However, if (\ref{rate}) holds
then in particular it follows that for $t_n=(n/c)^{-1/\alpha}$ we
have $\|Du(t_n)\|^2\ge n$. The fractal dimension of this sequence
$\{t_n\}$ is equal to $\alpha/(1+\alpha)$, and since we should have
$\alpha<1$, it follows that the dimension of this sequence should be
less than one half.

Rather than proving a strong estimate such as that in (\ref{rate}),
instead we prove that if we have a sequence $\|Du(t_n)\|^2\ge n$,
then the fractal dimension of $\{t_n\}$ can be no larger than one
half.

\begin{proposition}\label{tns}
  Suppose that
    $$
  \|Du(t_n)\|^2\ge n.
  $$
  for some sequence $t_n\ge 0$. Then $d_f(\{t_n\})\le 1/2$.
\end{proposition}

\def\T{{\mathscr T}}

In the proof we write ${\mathscr T}=\{t_n\}$.

\begin{proof}
  Suppose that $d_f(\T)=d>1/2$. Then for some $\delta$ with
  $1/2<\delta<d$ there exists a sequence $\epsilon_j\rightarrow0$
  such that
  $$
  N_j:=N({\mathscr T},\epsilon_j)>\epsilon_j^{-\delta},
  $$
  where $N(\T,\epsilon)$ denotes the maximal number of disjoint
  balls centred on $t_j$ of radius $\epsilon$. Let the centres of
  these balls be $t_{n_i}$, $1\le i\le N_j$. Since these balls are
  disjoint, we have
  $$
  \int_0^1\|Du(s)\|^2\,\d
  s\ge\sum_{i=1}^{N_j}\int_{t_{n_i}-\epsilon_j}^{t_{n_i}+\epsilon_j}\|Du(s)\|^2\,\d s>\sum_{i=1}^{N_j}\int_{t_{n_i}-\epsilon_j}^{t_{n_i}}\|Du(s)\|^2\,\d s.
  $$

Now, using the bound in (\ref{lower2}) it follows that
\bea
\int_{t-\epsilon}^t\|Du(s)\|^2\,\d s&\ge&\int_{t-\epsilon}^t
\frac{1}{\sqrt{\|Du(t)\|^{-4}+c(t-s)}}\,\d s\\
&=&\frac{1}{2c}\left[-\sqrt{\|Du(t)\|^{-4}+c(t-s)}\right]_{t-\epsilon}^t\\
&=&\frac{1}{2c}\left[\sqrt{\|Du(t)\|^{-4}+c\epsilon}-\|Du(t)\|^{-2}\right].
\eea

Noting that $\sqrt{X^2+c\epsilon}-X$ is a decreasing function of $X$
we therefore have
  $$
  \int_0^1\|Du(s)\|^2\,\d
  s\ge\frac{1}{2c}\,\sum_{i=1}^{N_j}\left[\sqrt{n_i^{-2}+c\epsilon_j}-n_i^{-1}\right],
  $$
and since
  $n_i\ge i$
  it follows that
  $$
  \int_0^1\|Du(s)\|^2\,\d s
  \ge\frac{1}{2c}\,\sum_{n=1}^{N_j}\left[\sqrt{n^{-2}+c\epsilon_j}-n^{-1}\right].
  $$

In order to proceed we note that
  $$
X<\sqrt{c\epsilon}/2\quad\Rightarrow\quad\sqrt{X^2+c\epsilon}-X>\frac{\sqrt{c\epsilon}}{2}.
$$
It follows, since $N_j\ge\epsilon_j^{-\delta}$ and $\delta>1/2$ we
have, for $j$ large enough,
\begin{eqnarray}
\int_0^1\|Du(s)\|^2\,\d s
  &\ge&\frac{1}{2c}\,\sum_{n=2(c\epsilon_j)^{-1/2}}^{\epsilon_j^{-\delta}}\left[\sqrt{n^{-2}+c\epsilon_j}-n^{-1}\right]\nonumber\\
  &\ge&\frac{1}{4\sqrt{c}}\,\sum_{n=2(c\epsilon_j)^{-1/2}}^{\epsilon_j^{-\delta}}\sqrt{\epsilon_j}\label{morecomplicated}\\
  &\ge&\frac{1}{4\sqrt{c}}[\epsilon_j^{(1/2)-\delta}-2/\sqrt{c}].\nonumber
  \end{eqnarray}
Since $\epsilon_j\rightarrow0$ as $j\rightarrow\infty$ and
$\delta>1/2$ it follows that $\int_0^1\|Du(s)\|^2\,\d s=+\infty$, a
contradiction.
\end{proof}

Our discussion above suggests that in fact one should expect to be
able to prove that $d_f(\T)<1/2$ provided one assumes in addition
that $t_n\downarrow0$ (an assumption not necessary for the proof of
the Proposition). This is further supported by the following
observation. If one has
\be{forn}
\|Du(n^{-1})\|^2 \geq n,
\ee
i.e.~$t_n=n^{-1}$ in the statement of Proposition \ref{tns}, then
although $d_f(\T)=1/2$ one can show that (\ref{forn}) is impossible.
Indeed, for all $s \in (t_{n+1}, t_n)$ we have
$$
\|Du(s)\|^2 \geq
\frac {\|Du(t_n)\|^2}{\sqrt{1+(t_n-s)\|Du(t_n)\|^4}} \geq
\frac{n}{\sqrt{1+\frac{n^2}{n(n+1)}}}.
$$
So
$$
\int_0^{t_1} \|Du(s)\|^2
\geq \sum _{n=1} ^{\infty}
\frac{1}{n(n+1)}\frac{n}{\sqrt{1+\frac{n^2}{n(n+1)}}} \geq \sum
_{n=1}^{\infty} \frac{1}{\sqrt{2}(n+1)} = + \infty,
$$
which is a
contradiction.

\section{The set of singular times}

A very similar, but slightly simpler argument, shows that if $\Sigma$ is the set of `singular times' of a weak solution,
i.e.~
$$
\Sigma=\{t\ge0:\ u(t)\notin V\}=\{t\ge0:\ \|Du(t)\|=\infty\},
$$
then $d_f(\Sigma)\le1/2$.

\begin{corollary}
  The upper box-counting dimension of the set $\Sigma$ of singular times satisfies $d_f(\Sigma)\le1/2$.
\end{corollary}

\begin{proof}
Suppose that $d_f(\Sigma)=d>1/2$. Then as in the proof of Proposition \ref{tns}, for some $\delta$ with $1/2<\delta<d$ there exists a sequence $\epsilon_j\rightarrow0$ such that $N_j=N(\Sigma,\epsilon_j)>\epsilon_j^{-\delta}$. If there centres of these balls are $t_n$, $1\le n\le N_j$, then since the balls are disjoint
  $$
  \int_0^1\|Du(s)\|^2\,\d
  s\ge\sum_{n=1}^{N_j}\int_{t_n-\epsilon_j}^{t_n+\epsilon_j}\|Du(s)\|^2\,\d s>\sum_{n=1}^{N_j}\int_{t_n-\epsilon_j}^{t_n}\|Du(s)\|^2\,\d s.
  $$
Since $\|Du(t_n)\|=+\infty$, the lower bound in (\ref{Lerays}) holds, and so
  $$
  \int_0^1\|Du(s)\|^2\,\d
  s\ge\frac{1}{2c}\,\sum_{i=1}^{N_j}\sqrt{c\epsilon_j}\ge\frac{1}{2\sqrt c}\,\epsilon_j^{(1/2)-\delta}.
  $$
  Since the right-hand side tends to infinity as $j\rightarrow\infty$, we must have $d_f(\Sigma)\le1/2$.
\end{proof}

Scheffer's well-known result that the $1/2$-dimensional Hausdorff
measure of the set of singular times is zero, ${\mathscr
H}^{1/2}(\Sigma)=0$, is based on the same ingredients as our result,
but neither is a consequence of the other. Indeed, in general the
box-counting dimension gives an upper bound for the Hausdorff
dimension, but from this observation the consequence of our result
is only that ${\mathscr H}^{1/2}(\Sigma)<+\infty$. However, note
that while any countable set $K$ has Hausdorff dimension zero (so in
particular ${\mathscr H}^s(K)=0$ for all $s>0$), the box-counting
dimension has the property that $d_f(X)=d_f(\overline{X})$, so it is easy to find a countable subset of the line with $d_f(X)=1$.
Thus our result $d_f(\Sigma)\le 1/2$ does indeed serve to limit
further the set of singular times.

\bigskip

\section*{Acknowledgments}
JCR is a Royal Society University Research Fellow, and would like to
thank the Society for all their support. WS is currently a visiting fellow in the Mathematics Institute at the University of Warwick under the Marie Curie Host Fellowship for the Transfer of Knowledge. We would both like to thank Jos\'e Rodrigo for some interesting and helpful conversations.

\References

\harvarditem{Cafarelli, Kohn, \& Nirenberg}{1982}{CKN} L.
Caffarelli, R. Kohn, \& L. Nirenberg. Partial regularity of suitable
weak solutions of the Navier-Stokes equations. {\it Comm. Pure Appl.
Math.}, 35: 771--831, 1982.

\harvarditem{Constantin \& Foias}{1988}{CF} P. Constantin and C.
Foias. {\it Navier-Stokes Equations}. University of Chicago Press,
Chicago, 1988.

\harvarditem{Doering \& Gibbon}{1995}{DG} C.R. Doering \& J.D.
Gibbon. {\it Applied Analysis of the Navier-Stokes Equations}.
Cambridge Texts in Applied Mathematics, Cambridge University Press,
Cambridge, 1995.

\harvarditem{Falconer}{1990}{Falc1990} K.J. Falconer. {\it Fractal
Geometry}. Wiley, Chichester, 1990.

\harvarditem{Langa \& Robinson}{2006}{LR} J.A. Langa \& J.C.
Robinson. Fractal dimension of a random invariant set. {\it J. Math.
Pures Appl.} 85: 269--294, 2006.

\harvarditem{Leray}{1934}{L34} J. Leray. Essai sur le mouvement d'un
fluide visqueux emplissant l'espace. {\it Acta Math.}, 63: 193--248,
1934.

\harvarditem{Robinson}{2001}{JCR} J.C. Robinson. {\it
Infinite-dimensional dynamical systems}. Cambridge Texts in Applied
Mathematics, Cambridge University Press, Cambridge, 2001.

\harvarditem{Robinson}{2006}{SEMA} J.C. Robinson. Regularity and
singularity in the three-dimensional Navier--Stokes equations. {\it
Bolet\'in de la Sociedad Espa\~nola de Matem\'atica Aplicada}, 35:
43--71.

\harvarditem{Scheffer}{1976}{Scheffer76} V. Scheffer. Turbulence and
Hausdorff dimension, in {\it Turbulence and Navier Stokes Equation,
Orsay 1975}, Springer LNM 565: 174--183, Springer Verlag, Berlin,
1976.

\harvarditem{Temam}{1977}{T77} R. Temam. {\it Navier-Stokes
Equations}. North Holland, Amsterdam, 1977. Reprinted by AMS
Chelsea, 2001.

\endrefs

\end{document}